\documentclass[a4paper,twoside,11pt]{article}

\usepackage{mypreamblepackage}

\begin{document}

\title{Surface Area and Curvature \\ of the General Ellipsoid} 
\author{D. Poelaert\protect\footnote{formerly: European Space Agency, ESTEC, Noordwijk,
The Netherlands.}\and J. Schniewind${}^*$\and F. Janssens${}^*$}\date{\copyright\ 2004}\maketitle

\section{Preliminaries}
%-------------------

Considering an ellipsoid whose principal semi-axes $a,b,c$ are ordered in such a way that \mb{$a\geqsl b \geqsl c \geqsl 0$},
we take the $(b,c)$-plane as the `equatorial' plane and axis $a$ as the `polar' axis.
A generic point on the surface of the ellipsoid is then represented in spherical coordinates $R,\Theta,\Phi$ where $R$ is the `radius'
(i.e. the length of the vector positioning the point from the centre), $\Theta$ is the `co-latitude' and $\Phi$ is the `longitude'.

Therefore, $\Theta$ is the central angle positioning
the generic point within the meridian ellipse passing through the point, whereas $\Phi$ is the central angle
positioning the equatorial projection of the generic point within the equatorial ellipse.
Obviously, $R$ will depend on the values chosen for $\Theta$ and $\Phi$. 

The central cartesian coordinates of the generic point are:
\beqna
X & = & R \cos\Theta \mb{\ ,}\nonumber\\
Y & = & R \sin\Theta \cos\Phi \mb{\ ,}\nonumber\\
Z & = & R \sin\Theta \sin\Phi \mb{\ .}\nonumber
\eeqna

When one uses the eccentric anomalies $\theta$ and $\varphi$ for the ellipses mentioned above, instead of their
central angle counterparts, the coordinates become:
\beqna
X & = & a \cos\theta \mb{\ ,}\nonumber\\
Y & = & b \sin\theta \cos\varphi\mb{\ ,} \nonumber\\
Z & = & c \sin\theta \sin\varphi \nonumber\mb{\ .}
\eeqna
The radius to a point upon the ellipsoid surface is then
\be
R=\sqrt{X^2\!+\!Y^2\!+\!Z^2}=\sqrt{a^2\cos^2\theta+b^2\sin^2\theta\cos^2\varphi+c^2\sin^2\theta\sin^2\varphi}\mb{\ .}
\ee{R}

One may verify the equation of the ellipsoid in cartesian coordinates:
\[\frac{X^2}{a^2}+\frac{Y^2}{b^2}+\frac{Z^2}{c^2}=1\mb{\ .}\]
Expressing the latter equation in terms of the central angles $\Theta$ and $\Phi$ provides another expression for
the radius at the generic point:
\be R=\frac{a b c}{\sqrt{b^2c^2\cos^2\Theta+c^2a^2\sin^2\Theta\cos^2\Phi+a^2b^2\sin^2\Theta\sin^2\Phi}}\mb{\ .}\ee{RR}

For each point upon the ellipsoid, there is another interesting quantity, namely the shortest distance $H$ (call it `height') from the centre
to the tangent plane to the ellipsoid at the point considered. Hence $H=R\cos\nu$, where $\nu$ is the acute angle between
the radius and the normal to the ellipsoid surface at the point. Differentiating the equation of the ellipsoid in cartesian coordinates, one 
obtains $\frac{X}{a^2}\dd X+\frac{Y}{b^2}\dd Y+\frac{Z}{c^2}\dd Z=0$. This shows that the vector with cartesian components
$\left[\frac{X}{a^2},\frac{Y}{b^2},\frac{Z}{c^2}\right]$ is orthogonal to the tangent plane and therefore normal to the surface;
let $N$ be its norm.
The equation for the ellipsoid also shows that the scalar product of this normal vector by the radius vector of cartesian
components $[X,Y,Z]$ must be equal to one; in other words, $N R \cos\nu=1$. Hence one must have
\beqna
H=\frac{1}{N} & = & \frac{1}{\sqrt{\frac{X^2}{a^4}+\frac{Y^2}{b^4}+\frac{Z^2}{c^4}}} \nonumber \\
                     & = &  \frac{a b c}{\sqrt{b^2c^2\cos{}^2\theta+c^2a^2\sin{}^2\theta\cos{}^2\varphi+
                               a^2b^2\sin{}^2\theta\sin{}^2\varphi}} \label{H} \\
                     & = & \frac{a b c\sqrt{b^2c^2\cos{}^2\Theta+c^2a^2\sin{}^2\Theta\cos{}^2\Phi+a^2b^2\sin{}^2\Theta\sin{}^2\Phi}}
                              {\sqrt{b^4c^4\cos{}^2\Theta+c^4a^4\sin{}^2\Theta\cos{}^2\Phi+a^4b^4\sin{}^2\Theta\sin{}^2\Phi}}\label{HH}
\mb{.}\eeqna 

\section{Surface Integral}\label{surface}
%-----------------------
Taking here the eccentric anomalies as variables, the general expression for the value of a surface integral is (see e.g. Ref.\cite{korn}
and Section \ref{curvature}):
\[  S=\int_0^{2\pi}\!\!\!\int_0^\pi\sqrt{U(\theta,\varphi) V(\theta,\varphi)-W^2(\theta,\varphi)}\dd\theta\dd\varphi \qquad\mb{where}\]
\beqnastar
U(\theta,\varphi)&=&\left[\left(\frac{\der X}{\der\theta}\right)^2+\left(\frac{\der Y}{\der\theta}\right)^2+
\left(\frac{\der Z}{\der\theta}\right)^2\right]\mb{\ ,}\\
V(\theta,\varphi)&=&\left[\left(\frac{\der X}{\der\varphi}\right)^2+\left(\frac{\der Y}{\der\varphi}\right)^2+
\left(\frac{\der Z}{\der\varphi}\right)^2\right]\mb{\ ,}\\
W(\theta,\varphi)&=&\left[\,\frac{\der X}{\der\theta}\frac{\der X}{\der\varphi}\,+\,\frac{\der Y}{\der\theta}\frac{\der Y}{\der\varphi}
\,+\,\frac{\der Z}{\der\theta}\frac{\der Z}{\der\varphi}\,\right]\mb{\ .}
\eeqnastar
Performing the calculation gives
\be
\!S\!=\!\!\int_0^{2\pi}\!\!\!\!\int_0^\pi\!\!\!\sqrt{b^2c^2\cos^2\theta\!+\!
c^2a^2\sin^2\theta\cos^2\varphi\!+\!a^2b^2\sin^2\theta\sin^2\varphi}\ \sin\theta\dd\theta\dd\varphi\mb{.} 
\ee{S}
One can perform a cyclic permutation of the semi-axes $a,b,c$ in Equ.(\ref{S}). This comes down to choosing
the $(c,a)$-plane or the $(a,b)$-plane as the equatorial plane; the angles $\theta$ and $\varphi$ belong of course to different planes in each
case. We found that the easiest way to calculate the surface integral necessitates the choice we have made; this, to avoid introducing 
complex numbers and, later on, elliptic functions with parameter greater than one or negative.

On the other hand, for expressing the surface integral, one can also use the central angles ($\Theta$,$\Phi$) and obtain, after calculation,
an alternate form to Equ.(\ref{S}), be it of more complicated nature:
\[
S=a^2b^2c^2\ *
\]
\be
    \!\int_0^{2\pi}\!\!\!\!\int_0^\pi\!\frac{\sqrt{b^4c^4\cos{}^2\Theta+c^4a^4\sin{}^2\Theta\cos{}^2\Phi+
    a^4b^4\sin{}^2\Theta\sin{}^2\Phi}}
    {\left[b^2c^2\cos{}^2\Theta+c^2a^2\sin{}^2\Theta\cos{}^2\Phi+a^2b^2\sin{}^2\Theta\sin{}^2\Phi\right]^2}
    \sin\Theta\dd\Theta\dd\Phi \mb{.}
\ee{SS}

\section{Other Theoretical relationships}
%-------------------------------------
Equ.(\ref{H}) shows immediately that Equ.(\ref{S}) for the area is in fact
\[
S=a b c \int_0^{2\pi}\!\!\!\int_0^\pi\frac{1}{H(\theta,\varphi)}\sin\theta\dd\theta\dd\varphi
\]
However, in Equ.(\ref{S}), $\theta$ and $\varphi$ are dummy variables; they may very well be \emph{interpreted} as being
$\Theta$ and $\Phi$, since the latter have the same intervals of integration as $\theta$ and $\varphi$!
By doing so, the surface integral may be written under the alternate form, after Equ.(\ref{RR}),
\[
S=a b c \int_0^{2\pi}\!\!\!\int_0^\pi \frac{1}{R(\Theta,\Phi)}\sin\Theta\dd\Theta\dd\Phi\mb{\ .}
\]
If one notices that the volume of the ellipsoid is $\Upsilon=4\pi a b c/3$ and that the infinitesimal element of solid angle $\dd\Psi$ in the
direction $(\Theta,\Phi)$ is, in spherical coordinates, $\sin\Theta\dd\Theta\dd\Phi$, one obtains the remarkable result that
\[
\frac{S}{\Upsilon}=3\int_0^{4\pi}\frac{1}{R(\Psi)}\frac{\dd\Psi}{4\pi}\mb{\ ,}
\]
where $\Psi$ is the value between $0$ and $4\pi$ of the solid angle positioning the generic point of radius $R$.

The equivalent at two dimensions of the latter expression is the ratio between the length $\Lambda$ of an ellipse and its area $A=\pi a b$,
using the central radius and the central angle:
\[
\frac{\Lambda}{A}=2\int_0^{2\pi}\frac{1}{R(\Phi)}\frac{\dd\Phi}{2\pi}\mb{\ .}
\]

Finally, consideration of Equs.(\ref{RR}),(\ref{HH}) and (\ref{SS}) shows that one can also write
\[ 
S=\int_0^{2\pi}\!\!\!\int_0^\pi\frac{R^3(\Theta,\Phi)}{H(\Theta,\Phi)}\sin\Theta\dd\Theta\dd\Phi
=\int_0^{4\pi}\frac{R^3(\Psi)}{H(\Psi)}\dd\Psi\mb{.}
\]

\section{Calculating the surface integral}
%-------------------------------------

Our starting integral will be Equ.(\ref{S}). If one notices that, by symmetry, the ellipsoid can be decomposed into eight similar 
\emph{octants}, we can limit the integration range for $\theta$ and $\varphi$ to $[0,\frac{\pi}{2}]$ and write
\[
S=8 a b\int_0^{\frac{\pi}{2}}\!\!\!\int_0^{\frac{\pi}{2}}\sqrt{\frac{c^2}{a^2}+\left[\left(\frac{c^2}{b^2}-\frac{c^2 }{a^2}\right)
+\left(1-\frac{c^2}{b^2}\right)\sin^2\varphi\right]\sin^2\theta}\ \sin\theta\dd\theta\dd\varphi\mb{\ .}
\]
We first let $\quad k(\varphi)=\sqrt{\left(\frac{c^2}{b^2}-\frac{c^2 }{a^2}\right)+\left(1-\frac{c^2}{b^2}\right)
\sin^2\varphi}\quad$ and $\quad e=\sqrt{1-\frac{c^2}{a^2}}\quad$, these quantities remaining within the interval $[0,1]$.
We obtain
\[
S=8 a b\int_0^\frac{\pi}{2}\!\!\!\int_0^\frac{\pi}{2}\sqrt{1-e^2+k^2(\varphi)\sin^2\theta}\ \sin\theta\dd\theta\dd\varphi\mb{\ .}
\]
In order to explicit the inner integral, we substitute $\xi=k\cos\theta$ and get
\[
\int_0^\frac{\pi}{2}\sqrt{1-e^2+k^2\sin^2\theta}\ \sin\theta\dd\theta=
\frac{1}{k}\int_0^k\sqrt{\left(1-e^2+k^2\right)-\xi^2}\dd\xi\mb{\ .}
\]
A primitive integral can be found in Ref.\cite{gr}, \S 2.262.1:
\beqnastar
\int\sqrt{\left(1-e^2+k^2\right)-\xi^2}\dd\xi & = & \frac{\xi\sqrt{\left(1-e^2+k^2\right)-\xi^2}}{2}+\\
   &  &  \frac{1-e^2+k^2}{2}\asin\frac{\xi}{\sqrt{\left(1-e^2+k^2\right)}}\mb{\ ,}
\eeqnastar
so that
\[
\int_0^\frac{\pi}{2}\sqrt{1-e^2+k^2\sin^2\theta}\ \sin\theta\dd\theta=\frac{\sqrt{1-e^2}}{2}+
\frac{1-e^2+k^2}{2k}\asin\frac{k}{\sqrt{1-e^2+k^2}}\mb{\ .}
\]
There are two occurrences of $k$ within the argument of the arcsine; for the range of values of $k$ and $e$ at hand,
we therefore better replace
\[
\asin\frac{k}{\sqrt{1-e^2+k^2}}=\atan\frac{k}{\sqrt{1-e^2}}\mb{\ .}
\]
The surface integral becomes thus
\[
S=2\pi b c+4 a b\int_0^\frac{\pi}{2}\frac{1-e^2+k^2(\varphi)}{k(\varphi)}\atan\frac{k(\varphi)}{\sqrt{1-e^2}}\dd\varphi
\]
or, in terms of variable $k$, using $k\dd k=(1-\frac{c^2}{b^2})\sin\varphi\cos\varphi\dd\varphi$,
\[
S=2\pi b c+4 a b\int_{\sqrt{\frac{c^2}{b^2}-\frac{c^2}{a^2}}}^{\sqrt{1-\frac{c^2}{a^2}}}
\frac{\left(1+\frac{a^2k^2}{c^2}\right)\atan\frac{a k}{c}}{\sqrt{\frac{a^2}{c^2}-\left(1+\frac{a^2k^2}{c^2}\right)}
\sqrt{\left(1+\frac{a^2k^2}{c^2}\right)-\frac{a^2}{b^2}}}\dd k\mb{\ .}
\]
In order to simplify this form, let $\eta=1+\left(\frac{a k}{c}\right)^2\geqsl 1$, so that
\be
S=2\pi b c+2 b c\int_\frac{a^2}{b^2}^\frac{a^2}{c^2}\frac{\eta\atan\sqrt{\eta -1}}{\sqrt{\eta -1}\sqrt{\frac{a^2}{c^2}-\eta}
\sqrt{\eta-\frac{a^2}{b^2}}}\dd\eta\mb{\ .}
\ee{Seta}
Expressing this integral in terms of known integrals will necessitate an integration by parts:
\[
\mb{let}\quad u(\eta)=\atan\sqrt{\eta-1} \quad \mb{and} \quad v'(\eta)=
\frac{\eta}{\sqrt{\eta-1}\sqrt{\frac{a^2}{c^2}-\eta}\sqrt{\eta-\frac{a^2}{b^2}}}\mb{\ .}
\]
First, we obtain easily $u'(\eta)=\frac{1}{2\eta\sqrt{\eta -1}}$.

Next, we need a primitive integral for $v'(\eta)$;
one could take $\int_{\frac{a^2}{b^2}}^\eta v'(\zeta)\dd\zeta$. However, by substracting from the latter the constant
$\int_{\frac{a^2}{b^2}}^{\frac{a^2}{c^2}}v'(\zeta)\dd\zeta$, one obtains another primitive integral, which is available in
Ref.\cite{gr} at \S 3.132.5: 

\be
v(\eta)=-\int_{\eta}^{\frac{a^2}{c^2}}v'(\zeta)\dd\zeta=-2\left[\rho^{-1}\textrm{F}(\tau(\eta),m)+
\rho\textrm{E}(\tau(\eta),m)\right]
\ee{v}

\[
\mb{with} \quad \rho=\sqrt{\frac{a^2}{c^2}-1} \quad \mb{,} \quad \tau(\eta)=\asin\left[\frac{b}{a}
\sqrt{\frac{a^2-c^2\eta}{b^2-c^2}}\right] \quad \mb{,} \quad m=\frac{a^2(b^2-c^2)}{b^2(a^2 -c^2)}
\]
and where $\textrm{F}(\tau,m)$ and $\textrm{E}(\tau,m)$ designate the elliptic integrals of the first and second kind, 
respectively. Quantities $\tau$ and $m$ always remain within the intervals $[0,\frac{\pi}{2}]$ and $[0,1]$, respectively.

Result (\ref{v}) can be checked directly by differentiation of $v(\eta)$. 
For that purpose, observe that
\[
v'=-2\left[\rho^{-1}\frac{\der\textrm{F}}{\der\tau}+\rho\frac{\der\textrm{E}}{\der\tau}\right]
\frac{\dd\tau}{\dd\eta}=-2\left[\frac{1}{\sqrt{\eta -1}}+\sqrt{\eta-1}\right]\frac{\dd\tau}{\dd\eta}
\]
\[
\mb{and that} \qquad \frac{\dd\tau}{\dd\eta}=-\frac{1}{2\sqrt{\frac{a^2}{c^2}-\eta}\sqrt{\eta-\frac{a^2}{b^2}}}\mb{\ .}
\]
With some skill, result (\ref{v}) can also be obtained by using the Mathematica\texttrademark\ software.

We may now express that
\[
\int_\frac{a^2}{b^2}^\frac{a^2}{c^2}u v'\dd\eta=\left[\rule{0pt}{3ex}u v\right]_\frac{a^2}{b^2}^\frac{a^2}{c^2}-
\int_\frac{a^2}{b^2}^\frac{a^2}{c^2}u' v\dd\eta\mb{\ .}
\]
One has $\ \tau(\frac{a^2}{c^2})=\asin 0=0\ $ and $\ \tau(\frac{a^2}{b^2})=\asin 1=\frac{\pi}{2}\ $.
Careful examination of the expression for $\tau$ using the definitions of $\eta$ and $k$ shows indeed that
$\ \tau=\asin\cos\varphi=\frac{\pi}{2}-\varphi\ $.
\[
\mb{Therefore} \quad \left[\rule{0pt}{2.9ex}u v\right]_\frac{a^2}{b^2}^\frac{a^2}{c^2}=
2\atan\sqrt{\frac{a^2}{b^2}-1}\left[\rho^{-1}\textrm{K}(m)+\rho\textrm{L}(m)\right]\mb{\ ,}
\]
where $\textrm{K}(m)$ and $\textrm{L}(m)$ are the complete elliptic integrals of the first and second kind, respectively.

Coming back to the surface integral, one has
\[
S = 2\pi b c+4 b c \atan\sqrt{\frac{a^2}{b^2}-1}\left[\rho^{-1}\textrm{K}(m)+\rho\textrm{L}(m)\right]+2 b c\; I
\]
where
\[
I=\int_\frac{a^2}{b^2}^\frac{a^2}{c^2}\frac{\rho^{-1}\textrm{F}\left(\tau(\eta),m\right)+
\rho\textrm{E}\left(\tau(\eta),m\right)}{\eta\sqrt{\eta-1}}\dd\eta\mb{\ .}
\]
With the knowledge that
\[
\eta=\frac{a^2}{c^2}-\left(\frac{a^2}{c^2}-\frac{a^2}{b^2}\right)\sin^2\tau \qquad \mb{and that}
\]
\[
\dd\eta=2\frac{a^2}{c^2}k\dd k=2(\frac{a^2}{c^2}-\frac{a^2}{b^2})\sin\varphi\cos\varphi\dd\varphi=
                                              -2(\frac{a^2}{c^2}-\frac{a^2}{b^2})\cos\tau\sin\tau\dd\tau\mb{\ ,}
\]

one can express integral $I$ in terms of $\tau$, after simplifications, as
\[
I=\frac{2}{\rho}\left(1-\frac{c^2}{b^2}\right)\int_0^\frac{\pi}{2}\frac{\left[\rho^{-1}\textrm{F}(\tau,m)+\rho\textrm{E}(\tau,m)\right]
\sin\tau\cos\tau}{(1-m e^2\sin^2\tau)\sqrt{1-m\sin^2\tau}}\dd\tau\mb{\ .}
\]

Fortunately, integral $I$ is now the combination of two definite integrals which can be found in Ref.\cite{gr}, the first one at \S 6.113.2,
\[
\int_0^\frac{\pi}{2}\frac{\textrm{F}(\tau,m)\sin\tau\cos\tau}{(1-m e^2\sin^2\tau)
\sqrt{1-m\sin^2\tau}}\dd\tau=
\]
\[-\frac{1}{m e\sqrt{1-e^2}}\left[\atan\left(\sqrt{1-m}\frac{e}{\sqrt{1-e^2}}\right)\textrm{K}(m)-\frac{\pi}{2}\textrm{F}(\asin e,m)\right]
\mb{,}\]
and the second one at \S 6.123,
\[
\int_0^\frac{\pi}{2}\frac{\textrm{E}(\tau,m)\sin\tau\cos\tau}{(1-m e^2\sin^2\tau)
\sqrt{1-m\sin^2\tau}}\dd\tau=
-\frac{1}{m e\sqrt{1-e^2}}\ \ *\]\[\left[\atan\left(\sqrt{1-m}\frac{e}{\sqrt{1-e^2}}\right)\!\textrm{L}(m)\!-\!
\frac{\pi}{2}\textrm{E}(\asin e,m)\!+\!\frac{\pi\sqrt{1-e^2}(1-\sqrt{1-m e^2})}{2 e}\right]\mb{.}
\]

Replacing all of these results into the surface integral and simplifying everything finally gives the surface area of the
ellipsoid:
\be
\boxed{S=2\pi\left[c^2+\frac{b c^2}{\sqrt{a^2-c^2}}\textrm{F}(\asin e,m)+b\sqrt{a^2-c^2}\textrm{E}(\asin e,m)\right]}
\ee{surf}
\[
\mb{with} \quad e=\sqrt{1-\frac{c^2}{a^2}} \quad \mb{and} \quad m=\frac{a^2(b^2-c^2)}{b^2(a^2-c^2)}\mb{\ .}
\]
This formula can be found without demonstration in the add-ons to Ref.\cite{math}, be it with a severe typographical error in versions
earlier than \mb{Mathematica\texttrademark\ 5}.

\section{Different kinds of ellipsoids}
%---------------------------------
Interpretation of the parameters $e$ and $m$ is important. Eccentricity $e$ is the largest eccentricity of all possible ellipses obtained
by cutting the ellipsoid by planes. If we let $m=\sin^2\gamma$, then $\gamma$ is the angle between the (b,c)-plane and any of the
two planes containing the \emph{circular} sections of the general ellipsoid. The circular sections are orthogonal to the (a,c)-ellipse.
\begin{quote}
When $\gamma<\atan\frac{a}{c}$, the ellipsoid shows a \emph{prolate} character.

When $\gamma>\atan\frac{a}{c}$, it shows an \emph{oblate} character.

When $\gamma=\atan\frac{a}{c}$,
the two planes containing the circular sections of the ellipsoid cut the (a,c)-ellipse along conjugated diameters (see Section \ref{curvature})
and the ellipsoid shows a \emph{spheroidal} character.  
\end{quote}

In order to derive all possible limit cases, let us rewrite formula (\ref{surf}) as a function of the three independent parameters
$a$, $e$ and $m$ only:
\be
\boxed{\!S\!\!=\!2\pi a^2\frac{\sqrt{1-e^2}}{\sqrt{1-m e^2}}\!\left[\!\sqrt{1-e^2}\sqrt{1-m e^2}+\!(1-e^2)\frac{\textrm{F}(\asin e,m)}{e}+
e\textrm{E}(\asin e,m)\right]\!}\mb{\,.}
\ee{newsurf} 

$\bullet$ For $m=0$ with $0<e<1$, one has $b=c$ and the ellipsoid is \emph{prolate of revolution} (about axis $a$). Since 
$\textrm{F}(\asin e,0)=\textrm{E}(\asin e,0)=\asin e$, the surface area (\ref{newsurf}) simplifies into
\be
\boxed{S_\textrm{prol.rev.}=2\pi a^2\sqrt{1-e^2}\left[\sqrt{1-e^2}+\frac{\asin e}{e}\right]}\mb{\ .}
\ee{prolrev}
One may verify that
\begin{eqnarray*}
S_\textrm{prol.rev.} & \xrightarrow[e\rightarrow 0]{} & 4\pi a^2 \qquad \mb{(sphere),} \\
S_\textrm{prol.rev.} & \xrightarrow[e\rightarrow 1]{} & 0 \qquad \mb{(bar-like shape along $a$).}
\end{eqnarray*}

$\bullet$ For $m=1$ with $0<e<1$, one has $a=b$ and the ellipsoid is \emph{oblate of revolution} (about axis $c$). With the knowledge that
$\int\frac{1}{\cos\varphi}\dd\varphi=\atanh\sin\varphi$, one has $\textrm{F}(\asin e,1)=\atanh e$ and $\textrm{E}(\asin e,1)=e$,
so that the surface area (\ref{newsurf}) simplifies into
\be 
\boxed{S_\textrm{obl.rev.}=2\pi a^2\left[1+(1-e^2)\frac{\atanh e}{e}\right]}\mb{\ .}
\ee{oblrev}
One may verify that
\[
S_\textrm{obl.rev.}\xrightarrow[e\rightarrow 0]\,4\pi a^2 \qquad \mb{(sphere).}
\]

$\bullet$ For both $m=1$ and $e=1$ though, the factor in front of Equ.(\ref{newsurf}) becomes indeterminate! But since, by definition,
\[
\frac{\sqrt{1-e^2}}{\sqrt{1-m e^2}}=\frac{b}{a}\mb{\ ,}
\]
this only means that $b$ may be chosen at will in the range $[0,a]$. Because $c=0$, we are indeed dealing here with a flat \emph{elliptic disc}
in the $(a,b)$-plane. Since $\lim_{e\rightarrow 1} \left[\left(1-e^2\right)\atanh e\right]=0$, its surface area is
\[
\boxed{S_\textrm{ell.disc} = 2\pi a b}\mb{\ ,}
\] 
the disc becoming circular when $b=a$.

$\bullet$ For $e=1$ together with $m$ indeterminate, one must have $c=b=0$, so that the ellipsoid degenerates into a \emph{bar} along
axis $a$ with, from (\ref{newsurf}),
\[
\boxed{S_\textrm{bar}=0}\mb{\ .}
\]

$\bullet$ For $e=0$, one has $c=b=a$ and the ellipsoid reduces to a \emph{sphere}. Parameter $m$ is indeterminate in this case but, since
$\lim_{e\rightarrow 0}\frac{\textrm{F}(\asin e,m)}{e}=1$ for all $m$, one obtains, from (\ref{newsurf}), as it should:
\[
\boxed{S_\textrm{sphere}=4\pi a^2}\mb{\ .}
\]

For completeness, we give in Table \ref{tab} a classification of all possible ellipsoids, together with their surface area $S$; we use the notations
\[
b^*=\sqrt{\frac{a^2+c^2}{2}} \qquad \mb{and} \qquad m^*=\frac{a^2}{a^2+c^2}
\]
for the intermediate semi-axis and the parameter of the spheroid.

\bt{!hbtp}{|c||c|c|c||c|} \hline
ellipsoid & semi-axes & $e$ & $m$ & area $S$ \\ \hline\hline
prolate of revolution & $0<c=b<a$ & $0<e<1$ & $m=0$ & Equ.(\ref{prolrev}) \\ 
general prolate & $0<c<b<b^*<a$ & $0<e<1$ & $0<m<m^*$ & Equ.(\ref{newsurf}) \\ 
general spheroid & $0<c<b=b^*<a$ & $0<e<1$ & $m=m^*$ & Equ.(\ref{newsurf}) \\
general oblate & $0<c<b^*<b<a$ & $0<e<1$ & $m^*<m<1$ & Equ.(\ref{newsurf}) \\
oblate of revolution & $0<c<b=a$ & $0<e<1$ & $m=1$ & Equ.(\ref{oblrev}) \\
elliptic disc & $0=c<b<a$ & $e=1$ & $m=1$ & $2\pi a b$ \\
circular disc & $0=c<b=a$ & $e=1$ & $m=1$ & $2\pi a^2$ \\
bar & $0=c=b<a$ & $e=1$ & $0\leqsl m\leqsl 1$ & $0$ \\
sphere & $0<c=b=a$ & $e=0$ & $0\leqsl m\leqsl 1$ & $4\pi a^2$ \\
point & $0=c=b=a$ & $0\leqsl e\leqsl 1$ & $0\leqsl m\leqsl 1$ & $0$ \\ \hline
\et{Classification of the ellipsoids}{tab}

\section{Curvature of the ellipsoid} \label{curvature}
%--------------------------------
A description of how to obtain the (two) curvatures at one point of a surface is given in Ref.\cite{korn}.
Using the position vector \bd{R} in coordinates $(\theta,\varphi)$, one first needs vectors $\bd{R}_\theta$ and $\bd{R}_\varphi$
which are the derivatives of vector \bd{R} with respect to $\theta$ and $\varphi$, respectively. One may then define the
(positive definite) matrix of metric coefficients 
\[ 
\ba{cc} U & W \\ W & V \ea = \ba{l} \bd{R}_\theta \\ \bd{R}_\varphi \ea \bd{\cdot} \ba{ll} \bd{R}_\theta & \bd{R}_\varphi \ea \mb{.}
\] 
\beqnastar
\mb{After calculation, one finds}\qquad U & = & a^2\sin^2\theta+\cos^2\theta\left(b^2\cos^2\varphi+c^2\sin^2\varphi\right) \mb{,}\\
V  & = & \sin^2\theta\left(b^2\sin^2\varphi+c^2\cos^2\varphi\right) \mb{,}\\
W & = & -\left(b^2-c^2\right)\sin\theta\cos\theta\sin\varphi\cos\varphi \mb{.}
\eeqnastar
Note that the surface area is $S=\int_0^{2\pi}\!\!\!\int_0^\pi\|\bd{R}_\theta\w\bd{R}_\varphi\|\dd\theta\dd\varphi$ (see Section \ref{surface}).
Next, knowing that the unit vector normal to the ellipsoid surface is
\[
\bd{n}=\frac{\bd{R}_\theta\w\bd{R}_\varphi}{\|\bd{R}_\theta\w\bd{R}_\varphi\|} \mb{\ ,}
\]
one needs its derivatives $\bd{n}_\theta$ and $\bd{n}_\varphi$ as well, in order to construct the (symmetric) matrix
\[
\ba{cc} \varkappa & \mu \\ \mu & \lambda \ea = \ba{l} \bd{R}_\theta \\ \bd{R}_\varphi \ea \bd{\cdot}
\ba{ll} \bd{n}_\theta & \bd{n}_\varphi \ea \mb{.}
\]
\beqnastar
\mb{One finds that}\quad \varkappa & = & \frac{a b c}
                                 {\sqrt{b^2c^2\cos^2\theta+c^2a^2\sin^2\theta\cos^2\varphi+a^2b^2\sin^2\theta\sin^2\varphi}}\mb{\ ,} \\
                                 \lambda & = & \varkappa\sin^2\theta \mb{\ ,}\\
                                 \mu       & = & 0 \mb{.}
\eeqnastar
One can now define the \emph{curvature tensor}
\be
\boxed{\bd{\chi}=\ba{ll} \bd{n}_\theta & \bd{n}_\varphi \ea \ba{cc} U & W \\ W & V \ea^{-1} \ba{l} \bd{R}_\theta \\ \bd{R}_\varphi \ea}
\ee{tensor}
which is a symmetric tensor because
\[
\ba{l}\bd{R}_\theta\\\bd{R}_\varphi\ea\bd{\cdot\chi\cdot}\ba{ll}\bd{R}_\theta&\bd{R}_\varphi\ea=\ba{cc}
\varkappa & \mu \\ \mu & \lambda\ea\mb{\ is a symmetric matrix.}
\]
From $\bd{n}_\theta=\bd{\chi\cdot R}_\theta$ and $\bd{n}_\varphi=\bd{\chi\cdot R}_\varphi$, one concludes that
\be
\boxed{\dd\bd{n}=\bd{\chi\cdot}\dd\bd{R}}\mb{\ .}
\ee{dn}
Hence, for a given $\dd\bd{R}$, the more pronounced the curvature, the more vector $\bd{n}$ changes its orientation!

The (scalar) curvature $\Chi{\ell}$ for the direction $\bd{\ell}=\frac{\dd\bd{R}}{\|\dd\bd{R}\|}$ is then
\[
\Chi{\ell}=\bd{\ell\cdot\chi\cdot\ell}=\frac{\dd\bd{R\cdot\chi\cdot}\dd\bd{R}}{\dd\bd{R}\bd{\cdot}\dd\bd{R}}=
\frac{\dd\bd{R}\bd{\cdot}\dd\bd{n}}{\dd\bd{R}\bd{\cdot}\dd\bd{R}}=
\frac{\ba{ll}\dd\theta&\dd\varphi\ea\ba{cc}\varkappa&\mu\\\mu&\lambda\ea\ba{l}\dd\theta\\ \dd\varphi\ea}
{\ba{ll}\dd\theta&\dd\varphi\ea\ba{cc}U&W\\W&V\ea\ba{l}\dd\theta\\ \dd\varphi\ea}\mb{.}
\] 
When $\dd\bd{R}$ is along a \emph{principal direction of curvature} (eigenvector of $\bd{\chi}$), then $\dd\bd{n}=\chi\dd\bd{R}$,
where the scalar $\chi$ is a \emph{principal curvature} (eigenvalue of $\bd{\chi}$). In such a case,
\[ 
\ba{l}\bd{R}_\theta\\\bd{R}_\varphi\ea\bd{\cdot}\ba{ll}\bd{n}_\theta&\bd{n}_\varphi\ea\ba{l}\dd\theta\\ \dd\varphi\ea=
\chi\ba{l}\bd{R}_\theta\\\bd{R}_\varphi\ea\bd{\cdot}\ba{ll}\bd{R}_\theta&\bd{R}_\varphi\ea\ba{l}\dd\theta\\ \dd\varphi\ea
\]
\[
\mb{or equivalently\quad} \left(\ba{cc}\varkappa & \mu\\ \mu & \lambda\ea-\chi\ba{cc}U & W\\W & V\ea\right)\ba{l}\dd\theta\\ \dd\varphi\ea=
\ba{c}0\\0\ea\mb{.}
\]
The characteristic equation of this eigenvalue problem allows one to calculate the two principal curvatures $\Chi{1}$ and $\Chi{2}$.
Their sum (twice the so-called \emph{mean curvature}) and their product (square of the so-called \emph{Gaussian curvature})
are found to be, using Equs.(\ref{R}) and (\ref{H}):
\be
\boxed{\Chi{1}+\Chi{2}=\frac{\varkappa V+\lambda U-2\mu W}{U V-W^2}=\frac{H^3\left(a^2+b^2+c^2-R^2\right)}{a^2b^2c^2}}
\ee{sum}
\be
\mb{and}\quad\boxed{\Chi{1}\Chi{2}=\frac{\varkappa\lambda-\mu^2}{U V-W^2}=\frac{H^4}{a^2b^2c^2}}\mb{\ ,}
\ee{product}
the latter expressions being independent of the choice of coordinates: $(\theta,\varphi)$ or $(\Theta,\Phi)$.
Being the eigenvectors of the symmetric tensor $\bd{\chi}$, the principal directions of curvature $\bd{\ell}_1$,$\bd{\ell}_2$ will
be orthonormal.

When, in a general ellipsoid, one translates a circular section parallel to itself, one obtains new circular sections whose centres
travel in the (a,c)-ellipse along a diameter which is said to be \emph{conjugated} to the diameter through which the circular section of
radius $b$ passes. When the plane of one of these parallel circular sections becomes tangent to the ellipsoid, the circular section becomes
a single point called an \emph{umbilic}. 
There are four such points on a general ellipsoid, with cartesian coordinates
\[
X=\pm a\sqrt{\frac{a^2-b^2}{a^2-c^2}} \quad \mb{,} \quad Y=0 \quad \mb{,} \quad Z=\pm c\sqrt{\frac{b^2-c^2}{a^2-c^2}}
\]
so that, for these umbilics,
\[
R=\sqrt{X^2+Y^2+Z^2}=\sqrt{a^2+c^2-b^2} \quad \mb{and} \quad H=\frac{1}{\sqrt{\frac{X^2}{a^4}+\frac{Y^2}{b^4}+
\frac{Z^2}{c^4}}}=\frac{a c}{b} \mb{\,.}
\]
With the use of Equs (\ref{sum}) and (\ref{product}), one can calculate the principal curvatures at the umbilics:
\[
\Chi{1}=\Chi{2}=\frac{a c}{b^3} \mb{\ .}
\] 
In accordance with the definition given in Ref.\cite{korn}, the curvature at an umbilical point is the same in any direction.

\end{document}